\documentclass{amsart}
\usepackage{amssymb}
\usepackage{eucal}
\usepackage{amsmath}
\usepackage{amscd}

\newtheorem{theorem}{Theorem}[section]
\newtheorem{proposition}[theorem]{Proposition}
\newtheorem{definition}[theorem]{Definition}
\newtheorem{lemma}[theorem]{Lemma}
\newtheorem{corollary}[theorem]{Corollary}
\newtheorem{prop-def}{Proposition-Definition}[section]
\newtheorem{claim}{Claim}[section]

\newcommand{\nc}{\newcommand}

\renewcommand{\Bbb}{\mathbb}
\renewcommand{\frak}{\mathfrak}

\newcommand{\efootnote}[1]{}

\renewcommand{\textbf}[1]{}

\nc{\dfootnote}[1]{{}}          
\nc{\ffootnote}[1]{\dfootnote{#1}}
\nc{\mfootnote}[1]{{}}        
\nc{\mlabel}[1]{\label{#1}}  
\nc{\mkeep}[1]{\marginpar{{\bf #1}}} 

\nc{\bin}[2]{ (_{\stackrel{\scs{#1}}{\scs{#2}}})}  
\nc{\binc}[2]{ \left (\!\! \begin{array}{c} \scs{#1}\\
    \scs{#2} \end{array}\!\! \right )}  
\nc{\bincc}[2]{  \left ( {\scs{#1} \atop
    \vspace{-1cm}\scs{#2}} \right )}  
\nc{\bs}{\bar{S}}
\nc{\la}{\longrightarrow}
\nc{\rar}{\rightarrow}
\nc{\dar}{\downarrow}
\nc{\dap}[1]{\downarrow \rlap{$\scriptstyle{#1}$}}
\nc{\uap}[1]{\uparrow \rlap{$\scriptstyle{#1}$}}
\nc{\defeq}{\stackrel{\rm def}{=}}
\nc{\disp}[1]{\displaystyle{#1}}
\nc{\dotcup}{\ \displaystyle{\bigcup^\bullet}\ }
\nc{\hcm}{\ \hat{,}\ }
\nc{\ot}{\otimes}
\nc{\hts}{\hat{\otimes}}
\nc{\hcirc}{\hat{\circ}}
\nc{\lleft}{[}
\nc{\lright}{]}
\nc{\curlyl}{\left \{ \begin{array}{c} {} \\ {} \end{array}
    \right .  \!\!\!\!\!\!\!}
\nc{\curlyr}{ \!\!\!\!\!\!\!
    \left . \begin{array}{c} {} \\ {} \end{array}
    \right \} }
\nc{\longmid}{\left | \begin{array}{c} {} \\ {} \end{array}
    \right . \!\!\!\!\!\!\!}
\nc{\ora}[1]{\stackrel{#1}{\rar}}
\nc{\ola}[1]{\stackrel{#1}{\la}}
\nc{\oeq}[1]{\stackrel{#1}{=}}
\nc{\scs}[1]{\scriptstyle{#1}}
\nc{\mrm}[1]{{\rm #1}}
\nc{\margin}[1]{\marginpar{\rm #1}}   
\nc{\dirlim}{\displaystyle{\lim_{\longrightarrow}}\,}
\nc{\invlim}{\displaystyle{\lim_{\longleftarrow}}\,}
\nc{\mvp}{\vspace{0.3cm}}
\nc{\tk}{^{(k)}}
\nc{\tp}{^\prime}
\nc{\ttp}{^{\prime\prime}}
\nc{\svp}{\vspace{2cm}}
\nc{\vp}{\vspace{8cm}}
\nc{\proofbegin}{\noindent{\bf Proof: }}
\nc{\proofend}{$\blacksquare$ \vspace{0.3cm}}
\nc{\modg}[1]{\!<\!\!{#1}\!\!>}
\nc{\intg}[1]{F_C(#1)}
\nc{\lmodg}{\!<\!\!}
\nc{\rmodg}{\!\!>\!}
\nc{\cpi}{\widehat{\Pi}}
\nc{\sha}{{\mbox{\cyr X}}}  
\nc{\shpr}{\diamond}    
\nc{\vep}{\varepsilon}
\nc{\labs}{\mid\!}
\nc{\rabs}{\!\mid}

\nc{\Hom}{\mrm{Hom}}
\nc{\id}{\mrm{id}}
\nc{\im}{\mrm{im}}
\nc{\supp}{\mathrm SSupp}

\nc{\Alg}{\mathbf{Alg}}
\nc{\Bax}{\mathbf{Bax}}
\nc{\bfk}{{\bf k}}
\nc{\bfone}{{\bf 1}}
\nc{\base}[1]{\bfone^{\otimes ({#1}+1)}} 
\nc{\Rings}{\mathbf{Rings}}
\nc{\Sets}{\mathbf{Sets}}

\nc{\BA}{{\Bbb A}}
\nc{\CC}{{\Bbb C}}
\nc{\DD}{{\Bbb D}}
\nc{\EE}{{\Bbb E}}
\nc{\FF}{{\Bbb F}}
\nc{\GG}{{\Bbb G}}
\nc{\HH}{{\Bbb H}}
\nc{\LL}{{\Bbb L}}
\nc{\NN}{{\Bbb N}}
\nc{\QQ}{{\Bbb Q}}
\nc{\RR}{{\Bbb R}}
\nc{\TT}{{\Bbb T}}
\nc{\VV}{{\Bbb V}}
\nc{\ZZ}{{\Bbb Z}}

\nc{\cala}{{\mathcal A}}
\nc{\calc}{{\mathcal C}}
\nc{\cald}{\mathcal{D}}
\nc{\cale}{{\mathcal E}}
\nc{\calf}{{\mathcal F}}
\nc{\calg}{{\mathcal G}}
\nc{\calh}{{\mathcal H}}
\nc{\cali}{{\mathcal I}}
\nc{\call}{{\mathcal L}}
\nc{\calm}{{\mathcal M}}
\nc{\caln}{{\mathcal N}}
\nc{\calo}{{\mathcal O}}
\nc{\calp}{{\mathcal P}}
\nc{\calr}{{\mathcal R}}
\nc{\calt}{{\mathcal T}}
\nc{\calw}{{\mathcal W}}
\nc{\calx}{{\mathcal X}}
\nc{\CA}{\mathcal{A}}

\nc{\fraka}{{\frak a}}
\nc{\frakA}{{\frak A}}
\nc{\frakB}{{\frak B}}
\nc{\frakm}{{\frak m}}
\nc{\frakp}{{\frak p}}

\font\cyr=wncyr10

\begin{document}

\title{Mixable Shuffles, Quasi-shuffles and Hopf Algebras
}
\thanks{MSC-class: 16A06, 05E99, 16W30}
\author{Kurusch Ebrahimi-Fard}
\address{Universit\"at Bonn -
         Physikalisches Institut,
         Nussallee 12,
         D-53115 Bonn,
         Germany}
\email{kurusch@ihes.fr}
\author{Li Guo}
\address{
Department of Mathematics and Computer Science,
Rutgers University,
Newark, NJ 07102, USA}
\email{liguo@newark.rutgers.edu}

\date{\today}


\begin{abstract}
The quasi-shuffle product and mixable shuffle product are both generalizations
of the shuffle product and have both been studied quite extensively recently.
We relate these two generalizations and realize quasi-shuffle
product algebras as subalgebras of mixable shuffle product algebras.
As an application, we obtain Hopf algebra structures in free Rota-Baxter
algebras.
\end{abstract}

\maketitle

\section{Introduction}
\label{intro}

This paper studies the relationship between the mixable shuffle product and
the quasi-shuffle product, both generalizations of the shuffle product.

Mixable shuffles arise from the study of Rota-Baxter algebras.
Let $\bfk$ be a commutative ring and let $\lambda \in \bfk$ be fixed.
A Rota-Baxter $\bfk$-algebra of weight $\lambda$ (previously called a Baxter
algebra)
is a pair $(R,P)$ in which $R$ is a $\bfk$-algebra and $P: R \to R$ is a
$\bfk$-linear map, such that
\begin{equation}
 P(x)P(y) = P(xP(y)) + P(P(x)y)+ \lambda P(xy),\ \forall x, y \in R.
\mlabel{eq:Ba}
\end{equation}

Rota-Baxter algebra was introduced by the mathematician Glen Baxter~\cite{Ba}
in 1960 to study the theory of fluctuations in probability.
Rota greatly contributed to the study of the Rota-Baxter
algebra by his pioneer work in the late 1960s and early 1970s~\cite{Ro1,Ro2,Ro3} and by
his survey articles in late 1990s~\cite{Ro4,Ro5}.
Unaware of these works, in the early 1980s the school around Faddeev, especially
Semenov-Tian-Shansky~\cite{STS},
developed a whole theory for the Lie algebraic version of equation (\ref{eq:Ba}),
which is nowadays well-know in the realm of the theory of integrable systems under the
name of (modified) classical Yang-Baxter equation.\footnote{The latter Baxter is the
Australian physicist Rodney Baxter.}
In recent years, Rota-Baxter algebras have found applications in
quantum field theory~\cite{C-K1,C-K2,E-G-K1,E-G-K2,E-G-K3},
dendriform algebras~\cite{A-L,EF,E-G1,Le}, number theory~\cite{Gu3},
Hopf algebras~\cite{A-G-K-O} and combinatorics~\cite{Gu2}.

Key to much of these
applications is the realization of the free objects in which the product is
defined by mixable shuffles~\cite{G-K1,G-K2} as a generalization of the shuffle
product.
The shuffle product is a natural generalization of the integration by parts
formula and its construction can be traced back to Chen's path integrals~\cite{Ch}
in 1950s.
It has been defined and studied in many areas of mathematics, such as Lie and Hopf
algebras, algebraic $K$-theory, algebraic topology and combinatorics. Its applications
can also be found in chemistry and biology. It naturally carries the notion of a
Rota-Baxter operator of weight zero.

Another paper on a generalization of the shuffle product was published~\cite{Ho2}
in the same year as the papers~\cite{G-K1,G-K2} on mixable shuffle products.
It was on the quasi-shuffle product by Hoffman\footnote{Hoffman
mentioned in~\cite{Ho2} that there was also a generalization in
the thesis of F. Fares~\cite{Fa}.}.
Hoffman's quasi-shuffle product plays a prominent r\^{o}le in
the recent studies of harmonic functions, quasi-symmetric functions,
multiple zeta values~\cite{Ho1,Ho3,H-O,B-B} (where in special cases
it is also called
stuffle product or harmonic product) and $q$-multiple zeta values~\cite{Br}.

Despite the extensive works on the two generalizations of shuffle products,
it appears that they were carried out without being aware of each other.
In particular,
the relation of quasi-shuffles with Rota-Baxter algebras seems unnoticed. For
example, in the numerous applications of quasi-shuffles in multiple zeta values
in the current literature, no connections with Rota-Baxter algebras and mixable shuffles
have been mentioned. In fact, concepts and results on Rota-Baxter algebras
were rediscovered in the study of multiple zeta values. For instance, the
construction of the stuffle product in~\cite{Br} follows easily from the
construction of free Rota-Baxter algebras in~\cite{Ca}, while the generalized
shuffle product in~\cite{Go} is the same as the mixable shuffle product
in~\cite{G-K1,G-K2}.

The situation is similar in the theory of dendriform algebras.
Even though both quasi-shuffles and Rota-Baxter algebras have been used to give examples
of dendriform algebras~\cite{A-L,LR}, no connection of the two have been made.
Also, in the work of
Kreimer, and Connes and Kreimer~\cite{Kr1,Kr2,C-K1,C-K2} on renormalization
theory in perturbative quantum field theory, both the shuffle and its
generalization in terms of the quasi-shuffle, and Rota-Baxter algebras
appeared, under different contexts.

It was noted in~\cite{EF1} that the two constructions should be related.
Our first goal of this paper is to make this connection precise.
We show that the recursive formula for the quasi-shuffle product has its
explicit form in the mixable shuffle product. Both can be derived from
the Baxter relation~(\ref{eq:Ba}) that defines a Rota-Baxter algebra
of weight 1.
We further show that the quasi-shuffle algebra on
a locally finite set is a subalgebra of a mixable shuffle algebra on the
corresponding locally finite algebra. With this connection, the concept of
quasi-shuffle algebras can be defined for a larger class of algebras.

This connection allows us to use the Hopf algebra structure on quasi-shuffle
algebras to obtain Hopf algebra structures on free Rota-Baxter algebras,
generalizing a previous work \cite{A-G-K-O} on this topic.
In the other direction, considering the critical r\^{o}le played by the quasi-shuffle
(stuffle) product in the recent studies of multiple zeta values and quasi-symmetric
functions, this connection should allow us to use the theory of Rota-Baxter algebras
in the studies of these exciting areas~\cite{E-G2}.

The paper is organized as follows. In the next section, we recall the concepts of
shuffles, quasi-shuffles and mixable shuffles, and describe their relations
(Theorem~\ref{thm:rel}).
In Section~\ref{sec:Hopf}, we use these connections to obtain Hopf algebra
structures on free Rota-Baxter algebras (Theorem~\ref{thm:main}).

\section{Shuffles, quasi-shuffles, and mixable shuffles}

For the convenience of the reader and for the ease of later references,
we recall the definition of each product before giving the relation among them.

\subsection{Shuffle product}
The shuffle product can be defined in two ways, one recursively, one
explicitly. We will see that Hoffman's quasi-shuffle product is a generalization
of the recursive definition and the mixable shuffle product is a generalization
of the explicit definition.

Let $\bfk$ be a commutative ring with identity $\bfone_\bfk$. Let $V$ be a $\bfk$-module. Consider the $\bfk$-module
\[
  T(V)=\bigoplus_{n\geq 0} V^{\otimes n}.
\]
Here the tensor products are taken over $\bfk$ and we take $V^{\otimes 0}=\bfk$.

Usually the shuffle product on $T(V)$ starts with the shuffles of permutations
\cite{Re,Sw}. For $m,n\in \NN_+$, define the set of {\bf $(m,n)$-shuffles} by
\[
  S(m,n)= \left \{ \sigma\in S_{m+n}
          \begin{array}{ll}
                {} \\ {}
          \end{array} \right .
\left |
\begin{array}{l}
 \sigma^{-1}(1)<\sigma^{-1}(2)<\ldots<\sigma^{-1}(m),\\
 \sigma^{-1}(m+1)<\sigma^{-1}(m+2)<\ldots<\sigma^{-1}(m+n)
\end{array}
\right \}.
\]
Here $S_{m+n}$ is the symmetric group on $m+n$ letters.

For $a = a_1 \otimes \ldots \otimes a_m \in V^{\otimes m}$,
$b = b_1 \otimes \ldots \otimes b_n \in V^{\otimes n}$ and
$\sigma\in S(m,n)$, the element
\[
  \sigma (a \otimes b) = u_{\sigma(1)}\otimes u_{\sigma(2)} \otimes
                         \ldots \otimes u_{\sigma(m+n)}\in V^{\otimes (m+n)},
\]
where
\[
    u_k = \left \{ \begin{array}{ll}
                   a_k,& \quad  1\leq k\leq m,\\
                   b_{k-m}, & \quad m+1\leq k\leq m+n, \end{array}
    \right .
\]
is called a {\bf shuffle} of $a$ and $b$. The sum
\begin{equation}
  a \:\sha\: b := \sum_{\sigma\in S(m,n)} \sigma(a\otimes b)
\mlabel{eq:shuf0}
\end{equation}
is called the {\bf shuffle product} of $a$ and $b$.
Also, by convention, $a \:\sha\: b$ is the scalar product if either $m=0$ or $n=0$.
The operation $\sha$ extends to a commutative and associative binary operation
on $T(V)$, making $T(V)$ into a commutative algebra with identity,
called the {\bf shuffle product algebra} generated by $V$.

The shuffle product on $T(V)$ can also be recursively defined as follows.
As above we choose two elements $a_1 \otimes \cdots \otimes a_m \in V^{\otimes m}$
and $b_1 \otimes \cdots \otimes b_n \in V^{\otimes n}$, and define
\begin{eqnarray*}
      a_0 \sha (b_1\otimes b_2 \otimes \ldots \otimes b_n) &=& a_0 b_1 \otimes b_2\otimes \ldots \otimes b_n, \\
     (a_1 \otimes a_2 \otimes \ldots \otimes a_m) \sha b_0 &=& b_0 a_1 \otimes a_2 \otimes \ldots \otimes a_m,
                                                                                   \ a_0,\, b_0 \in V^{\otimes 0}=\bfk,
\end{eqnarray*}
and
\begin{eqnarray}
      \lefteqn{(a_1\otimes  \ldots \otimes a_m)\sha (b_1 \otimes \ldots \otimes b_n) } \notag \\
                     &=& a_1 \otimes \big ((a_2\otimes \ldots \otimes a_m)\sha
             (b_1\otimes \ldots \otimes b_n)\big ) \mlabel{eq:rshuf1}\\
                     & & + b_1 \otimes \big ((a_1 \otimes \ldots \otimes a_m)\sha (b_2\otimes \ldots
                                              \otimes b_n)\big),\ a_i,\ b_j \in V. \notag
\end{eqnarray}

\begin{lemma}
For every element $v \in V$, the $\bfk$-linear map $P_{(v)}:(T(V),\sha) \to (T(V),\sha)$,
$P_{(v)}(a):=v \otimes a$ is a Rota-Baxter operator of weight zero.
\end{lemma}
\begin{proof}
    This is evident from the recursive definition of the shuffle product.
\end{proof}

\subsection{Quasi-shuffle product}

We recall the construction of quasi-shuffle algebras~\cite{Ho2}.
Let $X$ be a locally finite set, that is, $X$ is the disjoint union of
finite sets $X_n,\ n \geq 1$. The elements of $X_n$ are defined to have
degree $n$. Elements in $X$ are called letters and monomials in the letters
are called words. Even though the original paper~\cite{Ho2} only considered
$\bfk$ to be a subfield of $\CC$, much of the construction goes through for
any commutative ring $\bfk$. So we will work in this generality whenever
possible.
Consider the $\bfk$-module
underlying the noncommutative polynomial algebra $\frakA=\bfk\langle X\rangle$,
that is, the free $\bfk$-algebra generated by $X$. The identity 1 of $\frakA$
is called the empty word. Define $\bar{X}=X\cup \{0\}$. Suppose that there
is a pairing

\begin{equation}
   [\cdot,\cdot]:  \bar{X} \times \bar{X} \to \bar{X}
\mlabel{eq:set}
\end{equation}
with the properties
\begin{itemize}
 \item[S0.] $[a,0]=0$ for all $a\in \bar{X}$;
 \item[S1.] $[a,b]=[b,a]$ for all $a,b\in \bar{X}$;
 \item[S2.] $[[a,b],c]=[a,[b,c]]$ for all $a,b,c\in \bar{X}$;
 \item[S3.] either $[a,b]=0$ for all $a,b\in \bar{X}$,
or $\deg([a,b])=\deg(a)+\deg(b)$ for all $a,b\in \bar{X}$.
\end{itemize}
We define a {\bf Hoffman set} to be a locally finite set $X$ with a pairing
(\ref{eq:set}) that satisfies conditions S0-S3.

\begin{definition}
Let $\bfk$ be a commutative ring and let $X$ be a Hoffman set.
The {\bf quasi-shuffle product}
$*$ on $\frakA$ is defined recursively by
\begin{itemize}
 \item $1*w=w*1=w$ for any word $w$;
 \item $(aw_1) * (bw_2)=a(w_1*(bw_2))+b((aw_1)*w_2)+[a,b](w_1*w_2)$,
       for any words $w_1,w_2$ and letters $a,b$.
\end{itemize}
\mlabel{defn:quasi}
\end{definition}
When $[\cdot,\cdot]$ is identically zero, $*$ is the usual shuffle product
$\sha$ defined in Eq.~(\ref{eq:rshuf1}).

\begin{theorem} [(Hoffman)\cite{Ho2}]
\begin{enumerate}
 \item $(\frakA, *)$ is a commutative graded $\bfk$-algebra.
 \item When $[\cdot,\cdot]\equiv 0$, $(\frakA,*)$ is the shuffle product algebra
       $(T(V), \sha)$, where $V$ is the vector space generated by $X$.
 \item Suppose further $\bfk$ is subfield of $\CC$.
 Together with the coconcatenation comultiplication
\[
  \Delta: \frakA \to \frakA\otimes \frakA, w \mapsto \sum_{uv=w} u\otimes v
\]
where $uv$ is the concatenation of words, and counit
\[
\epsilon: \frakA \to \bfk, w\mapsto \delta_{w,1},
\]
$(\frakA,*)$ becomes a graded, connected bialgebra, in fact a Hopf algebra.
\end{enumerate}
\mlabel{thm:hoffman}
\end{theorem}

\subsection{Mixable shuffle product}

We next turn to the construction of mixable shuffle algebras and their
properties~\cite{G-K1}. The adjective mixable suggests that certain
elements in the shuffles can be mixed or merged.
We first give an explicit formula of the product before giving a recursive
definition which, under proper restrictions, will be seen to be equivalent
to Hoffman's quasi-shuffle product.

Intuitively, to form the shuffle product, one starts with two decks of cards
and puts together all possible shuffles of the two decks.
Suppose a shuffle of the two decks is taken and suppose a card from the first
deck is followed immediately by a card from the second deck, we allow the option
to merge the two cards and call the result a mixable shuffle.
To get the mixable shuffle product of the two decks of cards, one
puts together all possible mixable shuffles.

Given an $(m,n)$-shuffle $\sigma\in S(m,n)$, a pair of indices $(k, k+1)$,
$1\leq k< m+n$, is called an {\bf admissible pair} for $\sigma$ if
$\sigma(k)\leq m<\sigma(k+1)$. Denote $\calt^\sigma$ for the set of
admissible pairs for $\sigma$. For a subset $T$ of $\calt^\sigma$, call the
pair $(\sigma,T)$ a {\bf mixable $(m,n)$-shuffle}.
Let $\mid\! T\!\mid$ be the cardinality of $T$. By convention,
$(\sigma,T)=\sigma$ if $T=\emptyset$.
Denote
\begin{equation}
 \bs (m,n)=\{ (\sigma,T)\mid \sigma\in S(m,n),\
    T\subset \calt^\sigma\}
\mlabel{eq:mix}
\end{equation}
for the set of {\bf mixable $(m,n)$-shuffles}.

Let $A$ be a commutative $\bfk$-algebra not necessarily having an identity.
For $a = a_1 \otimes \ldots \otimes a_m \in A^{\otimes m}$,
$b = b_1 \otimes \ldots \otimes b_n \in A^{\otimes n}$ and $(\sigma,T)\in \bs(m,n)$,
the element
\[
  \sigma(a \otimes b; T)= u_{\sigma(1)}\hts u_{\sigma(2)} \hts
    \ldots \hts u_{\sigma(m+n)} \in A^{\otimes (m+n-\mid T\mid)},
\]
where for each pair $(k,k+1)$, $1\leq k< m+n$,
\[
   u_{\sigma(k)}\hts u_{\sigma(k+1)} =\left \{\begin{array}{ll}
          u_{\sigma(k)} u_{\sigma(k+1)},  &
                       \quad    (k,k+1)\in T\\
    u_{\sigma(k)}\otimes u_{\sigma(k+1)}, &
    \quad (k,k+1) \not\in T,
    \end{array} \right .
\]
is called a {\bf mixable shuffle} of the words $a$ and $b$.

Now fix $\lambda\in \bfk$. Define, for $a$ and $b$ as above,
the {\bf mixable shuffle product}
\begin{equation}
         a \shpr\!\!^+ b:\ = a \shpr\!\!^+_\lambda b\
                           = \sum_{(\sigma,T)\in \bs (m,n)} \lambda^{\mid T\mid } \sigma(a\otimes b;T)
                                                          \in \bigoplus_{k\leq m+n} A^{\otimes k}.
\mlabel{eq:shuf}
\end{equation}

As in the case of the shuffle product, the operation $\shpr^+ $ extends to a commutative
and associative binary operation on
\[
  \bigoplus_{k\geq 1} A^{\otimes k} = A\oplus A^{\otimes 2}\oplus \ldots
\]
Making it a commutative algebra without identity. Note that this is so even when $A$ has an
identity $\bfone_A$. In any case, we extend the product to
\begin{equation}
  \sha^+(A):=\sha_{\bfk,\lambda}^+(A):= \bigoplus_{k\in\NN} A^{\otimes k}
             =\bfk\oplus A\oplus A^{\otimes 2}\oplus \ldots,
\mlabel{eq:shaplus}
\end{equation}
making $\sha^+(A)$ a commutative algebra with identity $\bfone\in \bfk$~\cite{G-K1}.

Suppose $A$ has an identity $\bfone_A$.
Define
\begin{equation}
   \sha(A) := \sha_{\bfk,\lambda}(A) := A \otimes \sha_{\bfk,\lambda}^+(A)
\mlabel{eq:sha}
\end{equation}
to be the tensor product algebra, i.e., the {\bf augmented mixable shuffle product}
$\diamond$ on $\sha(A)$ is defined by:
\begin{equation}
(a_0 \otimes a) \diamond (b_0 \otimes b) := (a_0b_0) \otimes
                    (a \diamond^+ b),\ a_0,\ b_0\in A,\ a,\ b\in \sha^+(A).
\mlabel{eq:ashuf}
\end{equation}
Thus we have the algebra isomorphism (embedding of the second tensor factor)
\begin{equation}
\alpha: (\sha^+(A),\shpr^+) \to (\bfone_A\ot \sha^+(A),\shpr).
\mlabel{eq:double0}
\end{equation}

Define the $\bfk$-linear endomorphism $P_A$ on $\sha_\bfk(A)$ by assigning
\begin{eqnarray*}
  P_A( a_0\otimes a)
                &=&\bfone_A\otimes a_0\otimes a,\ a\in A^{\ot n},\ n\geq 1, \\
P_A(a_0\ot c)&=&\bfone_A \ot c a_0,\ c\in A^{\ot 0}=\bfk
\end{eqnarray*}
and
extending by additivity.
Let $j_A: A \rar \sha_\bfk(A)$ be the canonical inclusion map. Call $(\sha_\bfk(A),P_A)$
the {\bf {\rm (}mixable{\rm )} shuffle Rota-Baxter $\bfk$-algebra on $A$ of weight $\lambda$}.
The following theorem was proved in~\cite{G-K1}.

\begin{theorem}
\mlabel{thm:shua}
The shuffle Rota-Baxter algebra $(\sha_\bfk(A),P_A)$, together with the natural embedding
$j_A$, is a free Rota-Baxter $\bfk$-algebra on $A$ of weight $\lambda$. More precisely,
for any Rota-Baxter $\bfk$-algebra $(R,P)$ of weight $\lambda$ and algebra homomorphism
$f: A\to R$, there is a Rota-Baxter $\bfk$-algebra homomorphism
$\tilde{f}: (\sha_\bfk(A),P_A)\to (R,P)$ such that
$f=\tilde{f} \circ j_A$.
\end{theorem}
We will suppress $\bfk$ and $\lambda$ from $\sha_{\bfk,\lambda}(A)$ when there is no danger
of confusion.

\subsection{The connection}
\mlabel{ssec:conn}

We now establish the connection between quasi-shuffle product and mixable
shuffle product.

Let $\bfk$ be a commutative ring with identity.
Let $X=\cup_{n\geq 1} X_n$ be a Hoffman set.
Then the pairing $[\cdot,\cdot]$ in (\ref{eq:set}) extends by $\bfk$-linearity
to a binary operation on the free $\bfk$-module $A=\bfk\{X\}$ on $X$,
making $A$ into a commutative $\bfk$-algebra without
identity, with grading $A_n=\bfk\{X_n\}$, the free $\bfk$-module generated
by $X_n$.
Let $\tilde{A}=\bfk\oplus A$ be the unitary $\bfk$-algebra spanned by $A$.
Then $\tilde{A}=\bfk\{\tilde{X}\}$ where $\tilde{X}=\{ \bfone\} \cup X$
with $\bfone_{\tilde{A}}:=(\bfone_\bfk,0)$ the identity of $\tilde{A}$.
Here and in the rest of the paper, we will use $\bfone_A$ (instead of
$\bfone_{\tilde{A}}$) to denote this identity of $\tilde{A}$.
We will call $A$ (resp. $\tilde{A}$) the {\bf algebra
({\rm resp.} unitary algebra) generated by $X$}.

With the notations in Eq. (\ref{eq:shaplus}) and (\ref{eq:sha}), we have embeddings
\begin{equation}
\begin{array}{cccccc}
  \beta:& \sha^+(A) &\to & \sha^+(\tilde{A}) & \to & \sha(\tilde{A}),\, \\
         & a & \mapsto & a & \mapsto & \bfone_A \otimes a.
\end{array}
\mlabel{eq:embed}
\end{equation}
of $\bfk$-algebras.
Here the first embedding is induced by the embedding $A\hookrightarrow \tilde{A}$
and the second embedding is the natural one,
$\sha^+(\tilde{A}) \to \sha(\tilde{A}):= \tilde{A}\otimes \sha^+(\tilde{A}).$

\begin{theorem}
For a Hoffman set $X$,
the quasi-shuffle algebra $\frakA=\bfk\langle X\rangle$
is isomorphic to the algebra $\sha^+(A)$ and thus to the
subalgebra $\bfone_A\otimes \sha^+(A)$ of $\sha(\tilde{A})$ where the weight
$\lambda$ is 1. \mlabel{thm:rel}
\end{theorem}
\begin{proof}
We define
\[f: X\to X \subseteq A=A^{\ot 1} \subset \sha^+(A)\]
to be the canonical embedding.
We note that both $\frakA$, with the concatenation product, and
$\sha^+(A)$, with the tensor concatenation, are the free unitary non-commutative
$\bfk$-algebra on $X$. Thus $f$ extends uniquely to
a bijective map $\bar{f}: \frakA \to \sha^+(A)$
such that for any letters $a_1,\cdots, a_n\in X$, we have
\[\bar{f}(a_1\cdots a_n)=a_1\ot \cdots \ot a_n\in A^{\ot n}.\]
To distinguish elements, we will use $a=a_1\cdots a_n$ for an element in $\frakA$
and use $a^\ot:= a_1\ot \cdots \ot a_n$ for an element in $\sha^+(A)$ in the
rest of the proof.

To prove that $\bar{f}$ is also an isomorphism between $\frakA$, with the quasi-shuffle
product $\ast$, and $\sha^+(A)$, with the mixable shuffle product $\shpr^+$,
we just need to show that both products satisfy the same recursive relations.
We first note that the recursive relation of $\ast$ in Definition~\ref{defn:quasi}
can be rewritten as follows. For any $m,n\geq 1$ and
$a:=a_1\cdots a_m$, $b:=b_1\cdots b_n$ with
$a_i, b_j\in X,\ 1\leq i\leq m,\ 1\leq j\leq n,$ then
\\
\begin{equation}
\begin{array}{l}
 (1). \ a\ast b = a_1 b_1 + b_1 a_1 + [a_1,b_1], \ {\rm\ when\ } m,n=1,  \\
 (2). \ a \ast b = a_1  b_1  \cdots  b_n
 + b_1 \big(a_1\ast (b_2 \cdots b_n)\big) \\
\quad + [a_1, b_1]  b_2 \cdots  b_n, {\rm\ when\ } m=1, n\geq 2, \\
 (3).  \ a \ast b
= a_1 \big ((a_2 \cdots  a_m)\ast b_1 \big) + b_1 a_1 \cdots a_m  \\
\quad + [a_1,b_1]  a_2 \cdots  a_m,     {\rm\ when\ } m\geq 2, n=1, \\
 (4).  \ a \ast b = a_1  \big ((a_2 \cdots a_m)\ast
(b_1 \cdots  b_n)\big ) + b_1  \big ((a_1  \cdots  a_m)\ast (b_2 \cdots
 b_n)\big) \\
\quad   + [a_1, b_1]  \big ( (a_2 \cdots a_m) \ast
    (b_2 \cdots  b_n)\big ),\ {\rm\ when\ } m, n\geq 2.
\end{array}
\mlabel{eq:quasi2}
\end{equation}

On the other hand, consider the set $\bs (m,n)$ of mixable $(m,n)$-shuffles defined
in Eq. (\ref{eq:mix}). By \cite[Eq. (4)]{G-K1}, we have
\[
 \bs (m,n) = \bs_{1,0}(m,n) \dotcup
    \bs_{0,1}(m,n) \dotcup \bs_{1,1}(m,n)
\]
where
\begin{eqnarray*}
 \bs_{1,0}(m,n)&=& \{(\sigma,T)\in \bs(m,n) \mid (1,2)\not\in T,
    \sigma^{-1}(1)=1 \},\\
 \bs_{0,1}(m,n)&=& \{(\sigma,T)\in \bs(m,n)\mid (1,2)\not\in T,
    \sigma^{-1}(m+1)=1 \},\\
 \bs_{1,1}(m,n)&=&\{ (\sigma,T)\in \bs(m,n)\mid (1,2)\in T\}.
\end{eqnarray*}
Further, by \cite[p. 137]{G-K1}, for any $\lambda\in \bfk$,
\begin{eqnarray*}
\sum_{(\sigma,T)\in \bs_{1,0}(m,n)} \lambda^{|T|} \sigma (a^\ot \otimes b^\ot)
    &=&
a_1\otimes\!\!\! \sum_{(\sigma,T)\in \bs(m-1,n)} \lambda^{|T|}
    \sigma ((a_2\otimes \ldots \otimes a_m)\otimes b^\ot;T) \\
    &=&
\left \{ \begin{array}{ll}
a_1\otimes \big ((a_2\otimes \ldots \otimes a_m)\shpr^+ b^\ot \big),
& m\geq 2, \\
a_1\otimes b^\ot,
& m= 1. \end{array} \right .
\end{eqnarray*}
Similarly,
\begin{eqnarray*}
\sum_{(\sigma,T)\in \bs_{0,1}(m,n)} \lambda^{|T|} \sigma (a^\ot \otimes b^\ot)
    &=&
b_1\otimes\!\!\! \sum_{(\sigma,T)\in \bs(m,n-1)} \lambda^{|T|}
    \sigma (a^\ot \otimes (b_2\ot \cdots \ot b_n);T) \\
    &=&
\left \{ \begin{array}{ll}
b_1\otimes \big (a^\ot \shpr^+ (b_2\otimes \ldots \otimes b_n)\big),
& n\geq 2, \\
b_1 \ot a^\ot,
& n= 1; \end{array} \right .
\end{eqnarray*}
and
\begin{eqnarray*}
\lefteqn{ \sum_{(\sigma,T)\in \bs_{1,1}(m,n)} \lambda^{|T|}
\sigma (a^\ot \otimes b^\ot;T)}\\
&=& \sum_{(\sigma,T)\in \bs(m-1,n-1)}
  \lambda  [a_1, b_1] \otimes  \lambda^{|T|}\sigma((a_2\otimes \ldots \otimes a_{m})
    \otimes (b_2\otimes \ldots \otimes b_{n});T) \\
&=&
\left \{ \begin{array}{ll}
\lambda  [a_1,b_1] \otimes \big ((a_2\otimes \cdots \otimes a_{m}) \diamond^+
    (b_2\otimes \cdots \otimes b_{n} )\big ), &  m\geq 2, n\geq 2, \\
{\lambda [a_1, b_1]} \otimes (b_2\ot \cdots\ot  b_n), & m=1,n\geq 2, \\
{\lambda [a_1,b_1]}  \ot (a_2\ot \cdots \ot a_m), & m\geq 2, n=1,\\
{\lambda [a_1,b_1]},  & m=n=1. \end{array} \right .
\end{eqnarray*}
Since
\begin{eqnarray*}
\lefteqn{ a^\ot \shpr^+ b^\ot := \sum_{(\sigma,T)\in \bs(m,n)}
\lambda^{|T|}    \sigma (a^\ot\ot b^\ot)}\\
&=&
\sum_{(\sigma,T)\in \bs_{1,0}(m,n)} \lambda^{|T|} \sigma (a^\ot \ot b^\ot )+
\sum_{(\sigma,T)\in \bs_{0,1}(m,n)} \lambda^{|T|}  \sigma (a^\ot \ot b^\ot )\\
&& +
\sum_{(\sigma,T)\in \bs_{1,1}(m,n)} \lambda^{|T|} \sigma (a^\ot \ot b^\ot),
\end{eqnarray*}
we get
\begin{equation}
\begin{array}{l}
 (1). \ a^\ot \shpr^+ b^\ot = a_1 \ot b_1 + b_1 \ot a_1 + \lambda [a_1,b_1], \ {\rm\ when\ } m,n=1,  \\
 (2). \ a^\ot \shpr^+ b^\ot = a_1\ot  b_1 \ot  \cdots \ot b_n
 + b_1\ot \big(a_1\shpr^+ (b_2 \ot \cdots \ot b_n)\big) \\
\quad + \lambda [a_1, b_1]\ot  b_2 \ot \cdots\ot  b_n, {\rm\ when\ } m=1, n\geq 2, \\
 (3).  \ a^\ot \shpr^+ b^\ot
= a_1\ot \big ((a_2\ot \cdots \ot  a_m)\shpr^+ b_1 \big) + a_1 \ot \cdots \ot a_m \ot b_1 \\
\quad + \lambda [a_1,b_1]\ot  a_2 \ot \cdots \ot  a_m,     {\rm\ when\ } m\geq 2, n=1, \\
 (4).  \ a^\ot \shpr^+ b^\ot = a_1 \ot  \big ((a_2 \ot \cdots \ot a_m)\shpr^+
(b_1 \ot \cdots \ot  b_n)\big ) \\
\quad + b_1 \ot  \big ((a_1\ot  \cdots \ot  a_m)\shpr^+
(b_2 \ot \cdots \ot  b_n)\big) \\
\quad   + \lambda [a_1, b_1]\ot  \big ( (a_2 \ot \cdots \ot a_m) \shpr^+
    (b_2\ot \cdots \ot  b_n)\big ),\ {\rm\ when\ } m, n\geq 2.
\end{array}
\mlabel{eq:mixable}
\end{equation}
Thus when $\lambda=1$,
we have $\bar{f}(a\ast b) =a^\ot \shpr^+ b^\ot$ for all words $a$ and $b$
with $m,n\geq 1$, and hence for all $a$ and $b$ with $m,n\geq 0$ since when
$m=0$ or $n=0$, we have $a=1$ or $b=1$ and the multiplications through $\ast$ and
$\shpr^+$ are both given by the identity.
This proves the first isomorphism. The second one then follows from Eq. (\ref{eq:embed}).
\end{proof}
%
\begin{corollary}
Under the same assumptions of Theorem~\ref{thm:rel} and the additional
assumption that $\bfk$ is a subfield of $\CC$,
for any $\lambda\in \bfk$, the subalgebra $\sha^+(A)$ of $\sha^+(\tilde{A})$
and the subalgebra $\bfone_A\otimes \sha^+(A)$ of $\sha(\tilde{A})$ are Hopf
algebras.
\mlabel{co:hoffman}
\end{corollary}
In the next section, we will address the question on whether this Hopf algebra
can be extended to a larger Hopf algebra in $\sha(\tilde{A})$.
\begin{proof}
Because of the isomorphism (\ref{eq:double0}), we only need to prove the first
part of the statement.
When $\lambda=1$, this follows from Theorem~\ref{thm:rel} and Theorem~\ref{thm:hoffman}.
When $\lambda=0$, this is well-known (see \cite{Ho2}, for example).

Now assume $\lambda\neq 1,0$. We define a map
\[ g: \sha^+_\lambda (A) \to \sha^+_1(A)\]
by $g(a_1\ot \cdots \ot a_n)=\lambda^n a_1\ot \cdots \ot a_n, n\geq 1$
and $g(\bfone_\bfk)=\bfone_\bfk$.
Then by the first equation in (\ref{eq:mixable}), we have
\[ g(a_1\shpr^+_\lambda b_1)=g(a_1\ot b_1+b_1\ot a_1+\lambda [a_1,b_1])
=\lambda^2 (a_1\ot b_1 +b_1\ot a_1+[a_1,b_1])
=\lambda^2 (a_1\shpr^+_1 b_1)\]
which is just $g(a_1) \shpr^+_1 g(b_1). $
Using the other three equations in (\ref{eq:mixable}) and the induction on
$m+n$, we verify that
\[ g\big((a_1\ot \cdots \ot a_m)\shpr^+_\lambda (b_1\ot \cdots \ot b_n)\big)
=g(a_1\ot \cdots \ot a_m)\shpr^+_1 g(b_1\ot \cdots \ot b_n)\]
for all $m,n\geq 1$.
Thus we have $\sha^+_\lambda(A)$ is isomorphic to $\sha^+_1(A)$ and thus carries
a Hopf algebra structure.
\end{proof}

\section{Hopf algebras in Rota-Baxter algebras}
\mlabel{sec:Hopf}

We first recall the following theorem from~\cite{A-G-K-O}.
\begin{theorem} [Andrews-Guo-Keigher-Ono]
For any commutative ring $\bfk$ with identity and for any $\lambda \in \bfk$,
the free Rota-Baxter algebra $\sha_\lambda (\bfk)$ is a Hopf $\bfk$-algebra.
\mlabel{thm:alg}
\end{theorem}
As shown in~\cite{A-G-K-O}, when $\lambda=0$, we have the divided power Hopf algebra.

We now extend this result to $\sha(\tilde{A})$ for a $\bfk$-algebra $\tilde{A}$ coming
from a Hoffman set $X$.
To avoid confusion, we will use $\bfone_\bfk$ for the identity of $\bfk$
and $\bfone_A$ for the
identity of $\tilde{A}$ even though they are often identified under the
structure map $\bfk\to \tilde{A}$ of the unitary $\bfk$-algebra $\tilde{A}$.

Fix a $\lambda\in \bfk$. First note that, as a $\bfk$-module,
\[\sha^+(\bfk)= \bigoplus_{n\geq 0} \bfk^{\otimes n}=\bfk\oplus \bfk\oplus \bfk^{\ot 2} +\cdots.\]
There are two copies of $\bfk$ in the sum since $\bfk^{\ot 0}=\bfk=\bfk^{\ot 1}$. The identity of
$\sha^+(\bfk)$ is the identity in the first copy, which we  denote
by $\bfone_\bfk^{\ot 0}=\bfone$ as we did in (\ref{eq:shaplus}). Thus we have
\[ \sha^+(\bfk)= \bfk \bfone \oplus \bfk \bfone_\bfk \oplus
\bfk \bfone_\bfk^{\ot 2} \oplus \cdots
=\bigoplus_{n\geq 0} \bfk \bfone_\bfk^{\ot n}.\]
Then
\[\sha(\bfk) = \bfk\ot \sha^+(\bfk)= \bigoplus_{n\geq 0} \bfk(\bfone_\bfk \ot
\bfone_\bfk^{\ot n})\]
with the identity $\bfone_\bfk\ot \bfone$.
Since $\bfk$ is the base ring, the algebra homomorphism~(\ref{eq:embed}) gives
\[ \beta: (\sha^+(\bfk),\shpr^+) \cong (\sha(\bfk),\shpr).\]
Thus, by Theorem~\ref{thm:alg} we get
\begin{lemma}
For any $\lambda\in \bfk$, $(\sha^+(\bfk),\shpr^+)$ is a Hopf algebra.
\mlabel{lem:alg2}
\end{lemma}

For now let $\tilde{A}$ be any unitary $\bfk$-algebra with unit $\bfone_A$.
Then
\[ \sha^+(\tilde{A})= \bigoplus_{n\geq 0} \tilde{A}^{\ot n}
=\bfk\bfone \oplus \tilde{A}\oplus \tilde{A}^{\ot 2} \oplus \cdots\]
and
\[ \sha(\tilde{A})=\tilde{A}\ot \sha^+(\tilde{A})=
(\tilde{A}\ot \bfk\bfone) \oplus \tilde{A}^{\otimes 2} \oplus \tilde{A}^{\ot 3} \cdots. \]
Since $\sha(\tilde{A})$ is an $\tilde{A}$-algebra, and hence a $\bfk$-algebra, we have
the structure map $\gamma: \bfk\to \sha(\tilde{A})$ given by $\gamma(c)=c \bfone_A\ot \bfone$.
By the universal property of the free $\bfk$-Rota-Baxter algebra $\sha(\bfk)$, we have an
induced homomorphism $\gamma: \sha(\bfk)\to \sha(A)$ of Rota-Baxter algebras.
It is given by~\cite{G-K1}
\[ \gamma (\bfone_\bfk\ot \bfone_\bfk^{\ot n})
    =\bfone_A \ot \bfone_A^{\ot n},\ n\geq 0.\]
Let
\[\gamma^+: \sha^+(\bfk) \to \sha^+(\tilde{A}),\
\bfone_\bfk^{\ot n} \mapsto \bfone_A^{\ot n},\
n\geq 0.\]
We have the following commutative diagram
\[
\begin{array}{ccc}
\sha^+(\bfk) & \ola{\gamma^+} & \sha^+(\tilde{A}) \\
\downarrow & & \downarrow \\
\sha(\bfk)=\bfk\ot \sha^+(\bfk) & \ola{\gamma} & \sha(\tilde{A})=\tilde{A}\ot \sha^+(\tilde{A})
\end{array}
\]
where the vertical arrow are the injective maps to the second tensor factors.
\begin{theorem} 
Let $\bfk\subseteq \CC$ be a field. Let $X$ be a Hoffman set and let
the algebras $A$ and $\tilde{A}$ be the algebra and unitary algebra generated
by $X$ (as defined before Theorem~\ref{thm:rel}).
Let $\lambda\in \bfk$.
\begin{enumerate}
\item
The algebra product of
$\gamma^+(\sha^+(\bfk))$ and $\sha^+(A)$ in $\sha^+(\tilde{A})$ has a
Hopf algebra structure that expands the Hopf algebra structures on
$\gamma^+(\sha^+(\bfk))$ (see Lemma~\ref{lem:alg2}) and $\sha^+(A)$
(see Corollary~\ref{co:hoffman}).
\item
The algebra product of
$\gamma(\sha(\bfk))$ and $\bfone_A\ot\sha^+(A)$ in $\sha(\tilde{A})$ has a
Hopf algebra structure that expands the Hopf algebra structures on
$\gamma(\sha(\bfk))$ (see Theorem~\ref{thm:alg}) and $\bfone_A\ot \sha^+(A)$
(see Corollary~\ref{co:hoffman}).
\end{enumerate}
\mlabel{thm:main}
\end{theorem}
See Theorem~\ref{thm:main2} for a characterization of the elements in these
Hopf algebras.
\begin{proof}
Since the restriction of the isomorphism $\alpha: \sha^+(\tilde{A}) \to
\bfone_A\ot \sha^+(\tilde{A})$ in (\ref{eq:double0}) restricts to isomorphisms
$\gamma^+(\sha^+(\bfk)) \to \gamma(\sha(\bfk))$ and
$\sha^+(A)\to \bfone_A\ot \sha^+(A)$, we only need to prove the first statement.
Since the tensor product of two commutative, cocommutative Hopf algebras is
a Hopf algebra~\cite{Mo,L-N}, assuming Proposition~\ref{pp:tensor} which is stated
and proved below, we see
that $\gamma^+(\sha^+(\bfk))\shpr^+(\bfone\ot \sha^+(A))$
is a Hopf algebra for any $\lambda\in \bfk$ by Theorem~\ref{thm:alg} and
Corollary~\ref{co:hoffman}.
\end{proof}

\begin{proposition}
For any weight $\lambda\in \bfk$, let $\shpr^+$ be the mixable shuffle product
of weight $\lambda$. The two subalgebras $\gamma^+(\sha^+(\bfk))$ and
$\sha^+(A)$ of $\sha^+(\tilde{A})$ are linearly disjoint. Therefore,
$\gamma^+(\sha^+(\bfk))\shpr^+ \sha^+(A)$ is isomorphic to the tensor product
$\gamma(\sha(\bfk))\otimes (\bfone_A\ot \sha^+(A))$.
\mlabel{pp:tensor}
\end{proposition}

\begin{proof}
Let $\bfk\subseteq \CC, X, \tilde{X}, A$ and $\tilde{A}$ be as in
Theorem~\ref{thm:main}.
Since $X$ is locally finite, it is countable. So we can write
$X=\{y_n\, \big |\  n\geq 1\}$. Also denote $y_0=\bfone_A$, the unit of $\tilde{A}$.
Thus $\tilde{X}=\{y_n\,\big |\, n\in \NN\}$
and $\tilde{A}=\oplus_{n\geq 0} \bfk y_n.$
For $r\geq 1$ and $I=(i_1,\cdots,i_r)\in \NN^r$,
denote $y_{\otimes I}=y_{i_1}\otimes \cdots \otimes y_{i_r}$.
 Then
\[ A^{\ot r} =\bigoplus_{I\in \NN^r_{>0}} \bfk\, y_{\otimes I},
\quad \tilde{A}^{\ot r} =\bigoplus_{I\in \NN^r} \bfk\, y_{\ot I}.\]
By convention, we define $\NN^0=\NN^0_{>0}=\{\emptyset\}$,
and $y_{\otimes \emptyset}=\bfone$.
Let $\cali=\cup_{r\geq 0} \NN^r_{>0}$ and
$\tilde{\cali}=\cup_{r\geq 0} \NN^r$.
We then have
\[\sha^+(A)=\oplus_{I\in \cali}\, \bfk\, y_{\otimes I},
\quad
\sha^+(\tilde{A})=\oplus_{I\in \tilde{\cali}}\, \bfk\, y_{\otimes I}.\]

Recall that
\[\gamma^+(\sha^+(\bfk))=\oplus_{n\geq 0} \bfk \bfone_A^{\otimes n}.\]
So to prove that
$\oplus_{n\geq 0} \bfk\, \bfone_A^{\ot n}$ and $\sha^+(A)$
are linearly disjoint under the product $\shpr^+$,
we only need to prove
\begin{claim}
The set $\{ \bfone_A^{\otimes n}\shpr^+ y_{\otimes I} \,\big | \, n\geq 0,
I\in \cali\}$ is linearly independent.
\mlabel{cl:disjoint}
\end{claim}

Before proceeding further, we give a formula for the product
$\bfone_A^{\otimes n}\shpr^+ y_{\otimes I}$ which express a mixable shuffle product as
a sum of shuffle products in Eq.~(\ref{eq:shuf0}).
\begin{lemma}
For any $m\geq 0$ and $I\in \cali$, we have
\[
\bfone_A^{\ot m} \shpr^+ y_{\ot I} = \sum_{i=0}^{m} \lambda^i \binc{n}{i} \bfone_A^{\ot (m-i)}
\sha y_{\ot I}.
\]
\mlabel{lem:prod}
\end{lemma}
\begin{proof}
Define the length of $I\in \NN^r$ to be $\ell(I)=r$. We will prove by induction
on $w=m+\ell(I)$. When $w=0$, we have $m=\ell(I)=0$. Then
$\bfone_A^{\ot m}$ and $y_{\ot I}$ are both $\bfone$, so the lemma is clear,
as is if either $m=0$ or $\ell(I)=0$.
Suppose it holds for all $\bfone_A^{\ot m}\shpr^+ y_{\ot I}$ with $m+\ell<w$.
For given $\bfone_A^{\ot m}$ and $y_{\ot I}=y_{i_1}\ot \cdots \ot y_{i_r}$ with
$m\geq 1,\ r\geq 1$ and $m+r=w$, let $y'=y_{i_2}\ot \cdots \ot y_{i_r}$
if $r>1$ and $y'=\bfone$ if $r=1$.
Applying the recursive relation of $\shpr^+$ in Eq. (\ref{eq:mixable}), the induction
hypothesis, the Pascal equality and the recursive relation of $\sha$
in Eq.~(\ref{eq:rshuf1}), we have
\allowdisplaybreaks{
\begin{eqnarray*}
\lefteqn{\bfone_A^{\ot m} \shpr^+ y_{\ot I} = \bfone_A\otimes (\bfone_A^{\ot (m-1)} \shpr^+ y_{\ot I})
+ y_{i_1}\ot (\bfone_A^{\ot m} \shpr^+ y') + \lambda\, y_{i_1}\ot (\bfone_A^{\ot (m-1)} \shpr^+ y')}\\
&=& \bfone_A\ot \left ( \sum_{i=0}^{m-1} \lambda^i \binc{n}{i} \bfone_A^{\ot (m-1-i)}\sha y_{\ot I}
\right )
+ y_{i_1}\ot \left (\sum_{i=0}^{m} \lambda^i \binc{n-1}{i} \bfone_A^{\ot (m-i)}\sha y'
\right ) \\
&&+ \lambda\, y_{i_1}\ot \left (\sum_{i=0}^{m-1} \lambda^i \binc{n-1}{i} \bfone_A^{\ot (m-1-i)}\sha y'
\right ) \\
&=& \bfone_A\ot \left ( \sum_{i=0}^{m-1} \lambda^i \binc{n}{i} \bfone_A^{\ot (m-1-i)}\sha y_{\ot I}
\right )
+ y_{i_1}\ot \left (\sum_{i=0}^{m} \lambda^i \binc{n-1}{i} \bfone_A^{\ot (m-i)}\sha y'
\right ) \\
&&+ y_{i_1}\ot \left (\sum_{i=1}^{m} \lambda^i \binc{n-1}{i-1} \bfone_A^{\ot (m-i)}\sha y'
\right ) \\
&=& \bfone_A\ot \left ( \sum_{i=0}^{m-1} \lambda^i \binc{n}{i} \bfone_A^{\ot (m-1-i)}\sha y_{\ot I}
\right )
+ y_{i_1}\ot \left (\sum_{i=0}^{m} \lambda^i \binc{n}{i} \bfone_A^{\ot (m-i)}\sha y'
\right ) \\
&\oeq{(\ref{eq:rshuf1})}&
\sum_{i=0}^{m-1} \lambda^i \binc{n}{i} \bfone_A^{\ot (m-i)}
\sha y_{\ot I} + \lambda^m y\ot \binc{n}{m} \bfone\, \sha y'.
\end{eqnarray*}}
Since $y_{i_1}\ot (\bfone\, \sha y')=y_{i_1}\ot y' = y = \bfone_A^{\ot 0} \sha y$,
we get exactly what we want.
\end{proof}

We continue with the proof of Proposition~\ref{pp:tensor}.
For $r\geq 1$, let $[r]=(1,\cdots,r)$.
For a sequence $I=(i_1,\cdots,i_r)\in \NN^r$,
denote $\supp(I)$ (called sequential support) for the subsequence (with ordering)
of $I$ of non-zero entries.
For an all zero sequence $I=(0,\cdots,0)$ and the empty sequence $\emptyset$,
we define $\supp(I)=\emptyset$.
We then get a map
\[ \supp: \tilde{\cali} \to \tilde{\cali}.\]
Clearly, $\cali=\{I\in \tilde{\cali}\, \big |\, \supp(I)=I\}$.
So
\[ \tilde{\cali} = \dotcup_{I\in \cali} \supp^{-1} (I).\]
For each $I\in \cali$, consider the subset
$\calo_I=\{ y_{\ot J}\,\big |\, J\in \supp^{-1}(I)\}$.
Then we have
\[\{y_{\ot I} \, \big |\, I\in \tilde{\cali} \} = \dotcup_{I\in \cali} \calo_I.\]
So $\calo_I$ span
linearly independent subspaces of $\sha(\tilde{A})$.

By Lemma~\ref{lem:prod},
$\bfone_A^{\otimes n}\shpr^+ y_{\otimes I},\ n\geq 0,$ is in the linear
span of $\calo_I$.
Thus to prove Claim~\ref{cl:disjoint} and hence Proposition~\ref{pp:tensor},
we only need to prove that, for a fix $I\in \cali$, the subset
$\{ \bfone_A^{\otimes n}\shpr^+ y_{\otimes I}\,\big | \, n\geq 0\}$
is linearly independent.

Suppose the contrary.  Then there are integers $n_1>n_2>\cdots >n_r\geq 0$ and
$c_1\neq 0$ in $\bfk$ such that
$\sum_{i=1}^r c_i \bfone_A^{\otimes n_i}\shpr^+ y_I=0$.
Express this sum as a linear combination in terms of the basis $\calo_I$.
By Lemma~\ref{lem:prod}, the coefficient of
$\bfone_A^{\ot n_1}\otimes y_I$ is $c_1$, so we must have $c_1=0$, a contradiction.
\end{proof}

It is desirable to characterize the elements in the Hopf algebra
$\gamma^+(\sha^+(\bfk))\shpr^+\sha^+(A)$. This is our last goal in this article.
Recall that the {\bf length} of $y_{\ot I}$ with $I\in \NN^r,\ r\geq 0,$
is defined to be $\ell(y_{\ot r})=\ell(I)=r$.
For a given $I\in \NN^n$, the sum
$\sum y_{\otimes J}$ over $J\in \NN^n$ with $\supp(J)=
\supp (I)$, is called the {\bf one-shuffled element} of
$y_{\otimes I}$, denoted by $O(y_{\ot I})$. So $O(y_{\ot I})$ is the subset of
$\calo_I$ consisting of elements of length $\ell(I)$.
For example, if $I=(2,0,1)$, then the corresponding one-shuffle element of
$y_{\ot I}=y_2\ot \bfone_A \ot y_1$ is
$O(y_{\ot I})=y_2\ot \bfone_A \ot y_1 +\bfone_A \ot y_2\ot y_1 +y_2\ot y_1\ot \bfone_A.$
On the other hand, $O(y_{\ot I})$ is itself if $I$ is either an all zero sequence
or an all non-zero sequence.
It is so named because the sum can be obtained from shuffling the subsequence
of $y_{\ot I}$ of the $\bfone_A$-entries with the subsequence of $I$ of the
non-$\bfone_A$ entries (from $\supp(I)$).
To put it in another way, define a relation $\sim$ on $\tilde{\cali}$ by
$I_1\sim I_2$ if $\ell(I_1) =\ell(I_2)$ and $\supp  (I_1)=\supp (I_2)$. Then
it is easy to check that
$\sim$ is an equivalence relation and a one-shuffled element is of the form
$\sum y_{\ot J}$ where the sum is taken over all $J$ in an equivalence class.

We now give another version of Theorem~\ref{thm:main}.
\begin{theorem} 
Under the hypotheses of Theorem~\ref{thm:main},
the subspace of $\sha^+(\tilde{A})$ spanned by one-shuffled elements form
a Hopf algebra that  contains the Hopf algebras
$\gamma^+(\sha^+(\bfk))$ and $\sha^+(A)$.
\mlabel{thm:main2}
\end{theorem}

By Theorem~\ref{thm:main},
we only need to prove the following lemma.

\begin{lemma}
The product of $\gamma^+(\sha^+(\bfk))$ and $\sha^+(A)$ in $\sha^+(\tilde{A})$
is given by the subspace generated by one-shuffled elements.
\mlabel{pp:well-shuffle}
\end{lemma}

\begin{proof}
To prove the proposition, let $U$ be the product of $\gamma^+(\sha^+(\bfk))$ and
$\sha^+(A)$ in $\sha^+(\tilde{A})$,
and let $V$ be the subspace of one-shuffled elements
of $\sha^+(\tilde{A})$.
Then by Lemma~\ref{lem:prod} and the comments before the theorem, we have
$U\subseteq V$.
To prove $V\subseteq U$, we only need to show that, for each $k\geq 0$ and
$I\in (\NN^+)^n,\, n\geq 0$, the one-shuffled element
$\bfone_A^{\otimes k}\, \sha\, x_{\otimes I}$ is in $U$.
When $n=0$, $x_{\otimes I}=\bfone$. So
$\bfone_A^{\otimes k}\, \sha\, x_{\otimes I}=\bfone_A^{\otimes k}$ which is
in $\gamma^+(\sha^+(\bfk))$ and hence in $U$.
When $n\geq 1$, we use induction on $k$.
When $k=0$, then
$\bfone_A^{\otimes k}\, \sha\, x_{\otimes I} = x_{\otimes I}$ which is in
$\sha^+(A)$, hence is in $U$.
Assume that it is true for $\bfone_A^{\otimes k}, k< m$ and consider
$\bfone_A^{\otimes m}\, \sha\, x_{\otimes I}$.
By Lemma~\ref{lem:prod}, we have
\[
\bfone_A^{\ot m} \shpr^+ y_{\ot I} = \sum_{i=0}^{m} \lambda^i \binc{n}{i} \bfone_A^{\ot (m-i)}
\sha y_{\ot I}.
\]
The left hand side of the equation is in $U$ and, by induction, every term
on the right hand side except the first one (with $i=0$) is also in $U$.
Thus the first term, which is $\bfone_A^{\ot m}\sha y_{\ot}$, is also in $U$.
This completes the induction.
\end{proof}

\noindent
{\bf Acknowledgements. }
We thank Zongzhu Lin for helpful discussions. The first author thanks the I.H.\'E.S. for
the its warm hospitality, and the Ev. Studienwerk for financial support.

\end{document}